# BOX-KITES III:

# QUIZZICAL QUATERNIONS, MOCK OCTONIONS, AND OTHER ZERO-DIVISOR-SUPPRESSING "SLEEPER CELL" STRUCTURES IN THE SEDENIONS AND $2^N$-IONS

Robert P. C. de Marrais

**Abstract:** Building on two prior studies of zero-divisors (ZD's) generated by the Cayley-Dickson process, algebras we call "lariats" (**L**ine **A**lgebras of **R**eal and **I**maginary **A**xis **T**ransforms), linkable to quantum measurement, are discovered in the Sedenions, complementing the 7 isomorphic "box-kites" (pathway systems spanning octahedral lattices) interconnecting all primitive ZD's. By switching "edge signs," products among the diagonal line-pairs associated with each of a box-kite's 4 triangular, vertex-joined, "sails" generate not 6-cyclic ZD couplings when circuited, but 28 pairs of structures with Quaternionic multiplication tables – provided their symbols represent the *oriented diagonals as such*, not point-specifiable "units" residing on them. If a box-kite's 3 "struts" (pairs of opposite vertices, the *only* vertex pairings which do *not* contain mutual ZD's) each be combined with the ZD-free Quaternion copy uniquely associated with said box-kite, 21 lariats with *Octonionic* multiplication, one per each box-kite strut pair, are generated. Extending this approach to "emanation tables" (box-kite analogs in higher $2^n$-ions) indicates further ZD-masking "sleeper cell" structures, with "renormalization"'s basis possibly amenable to rethinking, thanks partly to the ZDs' newfound "Trip Sync" property, inhering throughout the $2^n$-ion hierarchy.

## 0. Zero-Divisors and Box-Kites Thus Far . . .

It's been common knowledge since the 19th Century that, by use of a simple algorithm dubbed the "Cayley-Dickson Process" (CDP), one can start with the real numbers, then progressively double dimensions to yield, in sequence, the Imaginaries, Quaternions, Octonions, Sedenions, and general $2^n$-ions. For almost as long, it's also been known that each such doubling causes increased flexibility in some senses, but loss of crucial structure in others: order relations and uniqueness of powers break down with Imaginaries (without whose *angular* power law, however, and remarkably simple relationship to stereographic projection, modern mathematics and wave mechanics would be unthinkable); Quaternions don't commute (thanks to which "failure," the vector-mechanical description of real forces in the physical world is derivable); Octonions are non-associative (making cause-effect relationships untraceable, but quark confinement and M-Theory's "brane scan" thinkable). Beyond these 8-D imaginaries, though, Hurwitz showed a century ago (see Kantor and Solodovnikov) that the very notions of "field" and a metrical "norm" fall apart, thanks to the appearance of a new kind of "monster": the *zero-divisor*.

Non-zero numbers whose product is zero at first seemed pathological, but are now – thanks to quantum mechanical "projection operators" and their formal link to "observables" – quite acceptable in "mixed" company. That is, if one defines a space with not just reals and imaginaries, but Pauli-spin-matrix-based "mirror numbers" too, whose non-real units square to +1, then the Clifford Algebras which allow them all to intercommunicate necessarily contain zero divisors. These, however, are of a special sort: if the "mirror number" unit be written as "e," then the two quantities ½ (1 ± e) are not just *idempotents* (hence algebraic analogs of "1"'s among Reals) on two orthogonal axes in a 4-D "bicomplex" space; more, they have *zero* for their product – which product obtains for *any* pair of points taken one from each axis. (As we'll see, the association of the zero-divisor property with *lines* of numbers, not just points, has quite far-reaching consequences.)

In "hypercomplex" spaces spanned only by the real axis plus lines whose units square to –1, however, the emergence of zero divisors is more problematic, and long remained poorly understood. In the 16-D Sedenions, the first extension of the imaginaries where they announce themselves, it took until 1998 before Guillermo Moreno showed a *homo*morphism between zero-divisor structure and the *auto*morphism group of the Octonions, whose Coxeter-Dynkin diagram is known in the argot of the "A, D, E" classification as "$G_2$" – and not just in the Sedenions, but in *all higher $2^n$-ions* (since $G_2$ is also the "derivation algebra" upon whose backbone all higher CDP imaginaries articulate themselves). But: homo- is not iso-, and concrete methods for extracting and studying zero-divisors were still lacking . . . until the Fall of 2000.

In the first paper in this series, a remarkably simple method was demonstrated for generating, and studying the interconnections among, *all* Sedenion zero-divisors. By employing (among the literally hundreds of billions of options) a scheme for labeling the 15 imaginary axes based strictly on XOR calculations (the index of $i_A$ times $i_B$ will always be A *xor* B), and conforming to CDP-generated signing patterns, zero divisors were shown to fill 84 diagonal lines, spanning 42 planes dubbed "Assessors," the axes of the latter being simply derived as follows: pick any arbitrary Octonion **'o'** (index 1-7, for 7 choices), and any pure Sedenion **'S'** whose index is greater than 8 and not equal to the XOR of **'o'** *with* **8** (for 6 choices); then all points on the orthogonal diagonals **o ± S** will be zero-divisors, but *never* with points on each other. Rather, diagonals from *different* Assessors will mutually zero-divide each other, according to easy-to-visualize path-selection rules.

The visualizing strategy was based on placing labels for each of the 42 Assessors on exactly *1* of the 6 vertices of *1* of 7 isomorphic octahedral lattices (the "Box-Kites"), then taking products of points along one or the other diagonal associated with each neighboring vertex, according to an "edge-signing"protocol. If the sign along the edge adjoining 2 vertices is "+," this indicates the product of *similarly* oriented diagonals (both either **o + S** or **o – S**: hence, "/· /" or "\· \"); if the sign is "–," then only pairings of *oppositely* oriented diagonals ("/ · \" or "\ · /") will mutually zero-divide. But how allocate the edge-signs? And, even more fundamentally, how partition the 42 Assessors among the 7 Box-Kites? Upon answering the second, we will see the first question almost answers itself.

We state without proof the *Strut Constant Rule*: each of the 7 ways of excluding 1 of the 7 Octonions, whose index we'll call the "strut constant," $s$, results in a set of 6 which can be mapped directly to 6 vertices of precisely one Box-Kite; moreover, the index of the Sedenion making an Assessor with each, is the XOR of the given Octonion with the "excess" $X = 8 + s$. Designate each Box-Kite with the Roman numeral of its strut constant. *Box-Kite I* thus has these 6 Assessor index-pairs at its vertices: $(2 \pm 11)$, $(3 \pm 10)$, $(4 \pm 13)$, $(5 \pm 12)$, $(6 \pm 15)$, $(7 \pm 14)$. (That is, 2 *xor* 11 = 3 *xor* 10 = 4 *xor* 13 = 5 *xor* 12 = 6 *xor* 15 = 7 *xor* 14 = $X$ = 8+1 = 9.) These 6 are further partitioned into mutually non-zero-dividing pairs (the *struts*): those 3 couplings whose Octonions' XOR would yield the missing strut-constant-indexed Octonion. With $s = 1$, the Assessor pairs with Octonions indexed respectively by 2 and 3, 4 and 5, and 7 and 6, would sit at opposite ends of diagonals (each playing that role for 2 orthogonal squares) on the octahedral lattice. We write 7 and 6 out of *counting* order, because that is their *signing order* in the most common XOR-based indexing scheme: the 7 Octonion triplets ("O-trips"), each forming a Quaternion copy with the reals, are listed, lowest index first, in left-right, then cyclical sequence, so that each successive pair's product in an O-trip is *positively* signed:

(1, 2, 3); (1, 4, 5); (1, 7, 6); (2, 4, 6); (2, 5, 7); (3, 4, 7); (3, 6, 5)

The strut-pairings just considered, then, are written in cyclic order as (2, 3); (4, 5); (7,6). The 3 terminal indices can always be rearranged to form an O-trip in *ascending signing order* (ASO) per the above list: (3, 6, 5) for Box-Kite I. Use this fact to mutually orient strut-pairs as follows: pick a triangle on the octahedron, label vertices A, B, C, and place "–" signs on the edge-lines between them. Label vertices at opposite ends of their struts F, E, D respectively, putting "–" signs on *their* triangle's edge-lines, and "+" signs on the 6 remaining edges. For Box-Kite I, set (A, B, C) to (3, 6, 5), hence (4, 7, 2) to F, E, D, in that order. O-trips now fill the vertices of 4 of the Box-Kite's 8 triangles, which we call "Sails." These share no edges, touching only at vertices, separated checkerboard-style by 4 "empty" triangles we call "Vents." The other Sails are (A, D, E), (F, D, B) and (F, C, E), these 3 having only 1 "–" edge-sign each. (All such letter triplets represent indices in signing order sequence – which, for Box-Kites I, II, III, are all "ASO" as well.)

Multiplying Assessors in sequence around a Sail's vertices, according to the edge-signing rule, a 6-cyclic listing of zero-divisor couplings results. Start with the $(- -, ++)$ diagonal at A; circuiting (A, B, C) in Box-Kite I yields this progression: $(i_3 + i_{10})(i_6 - i_{15})$ = $(i_6 - i_{15})(i_5 + i_{12})$ = $(i_5 + i_{12})(i_3 - i_{10})$ = $(i_3 - i_{10})(i_6 + i_{15})$ = $(i_6 + i_{15})(i_5 - i_{12})$ = $(i_5 - i_{12})(i_3 + i_{10})$ = 0. Observe, too, that in any Sail, each vertex's Octonion forms a Quaternionic axis-trio with the other 2 vertices' Sedenions, all 21 such "S-trips" with no "8" in them appearing *only* in such contexts, precisely once each among the 7 (A, B, C) Sails.

## 1. Groups, Loops and Lariats

For the reader's convenience, the 28 S-trips are listed in signing order below, followed by a "Strut Table" which shows Box-Kite structures per the above specifications. The only column in the table which will seem strange to those who haven't read the first monograph in this series (henceforth, "BK1," with its successor shorthanded "BK2") is the second: explaining what GoTo numbers are, though, will take us to our next topic, where new material will begin to reveal itself.

|          |           |           |           |
|----------|-----------|-----------|-----------|
| (1, 8, 9)  | (1, 11, 10) | (1, 13, 12) | (1, 14, 15) |
| (2, 8, 10) | (2, 9, 11)  | (2, 14, 12) | (2, 15, 13) |
| (3, 8, 11) | (3, 10, 9)  | (3, 15, 12) | (3, 13, 14) |
| (4, 8, 12) | (4, 9, 13)  | (4, 10, 14) | (4, 11, 15) |
| (5, 8, 13) | (5, 12, 9)  | (5, 10, 15) | (5, 14, 11) |
| (6, 8, 14) | (6, 15, 9)  | (6, 12, 10) | (6, 11, 13) |
| (7, 8, 15) | (7, 9, 14)  | (7, 13, 10) | (7, 12, 11) |

Note that the first column of S-trips contains the index-8 imaginary, and hence none of these appear in any Box-Kite Sails. They do, however, play a central role, since each contains precisely and only those indices forbidden from appearing in precisely 1 Box-Kite. In the first monograph, the offending index number's unique tagging of a "tame" kind of Quaternionic 3-sphere, in conjunction with the well-known billiards-based phrase for the unplayability of shots one is forced to take behind it, suggested calling such a zero-divisor-bereft S-trip an "8-ball." When its actions are treated in conjunction with those of its associated Box-Kite, we will see some surprises will be waiting for us. Before we get to those, though, here's the promised aid to Box-Kite building, followed with an explanation of the "GoTo Numbers" column.

### STRUT TABLE

| Box-Kite | GoTo Numbers | A | B | C | D | E | F |
|----------|--------------|-------|--------|-------|-------|--------|-------|
| I    | 7, 6, 4, 5 | **3**, 10 | **6**, 15 | **5**, 12 | 4, 13 | 7, 14 | 2, 11 |
| II   | 3, 2, 6, 7 | **1**, 11 | **7**, 13 | **6**, 12 | 4, 14 | 5, 15 | 3, 9  |
| III  | 5, 4, 2, 3 | **2**, 9  | **5**, 14 | **7**, 12 | 4, 15 | 6, 13 | 1, 10 |
| IV   | 1, 3, 5, 7 | **1**, 13 | **2**, 14 | **3**, 15 | 7, 11 | 6, 10 | 5, 9  |
| V    | 4, 1, 6, 3 | **2**, 15 | **4**, 9  | **6**, 11 | 3, 14 | 1, 12 | 7, 10 |
| VI   | 6, 1, 2, 5 | **3**, 13 | **4**, 10 | **7**, 9  | 1, 15 | 2, 12 | 5, 11 |
| VII  | 2, 1, 4, 7 | **1**, 14 | **4**, 11 | **5**, 10 | 2, 13 | 3, 12 | 6, 9  |

In the listing of O-trips in ASO sequence given above, the 7 were put in ascending order, as sorted by first digit, then second if necessary. The "GoTo numbers" refer to O-trips by this indexing scheme, indicating thereby 7 "Octonion counterfeits" whose non-real axes consist of the given O-trip plus the 4 Sedenions which aren't indexed by 8, the strut constant, or the XOR of these two. For each Box-Kite, 4 such numbers are listed, indicating in the order written the 4 "Octonion counterfeits" from which the 4 successive Sails (which we now list with edge-signs inserted) derive. For Box-Kite I, "7, 6, 4, 5" means the pattern of zero-divisors in A–B–C– is uniquely contained in the 7-dimensional space spanned by imaginaries with indices 3, 6, 5 (ABC's O-trip), plus 9, 10, 12 and 15; Sail A+D–E+ is embedded in the space whose axial indices are (3, 4, 7, 9, 10, 13, 14); Sail F–D+B+ is linked to (2, 4, 6, 9, 11, 13, 15); and F+C+E–, to (2, 5, 7, 9, 11, 12, 14).

In BK1, these 8-D algebras were called "automorphemes," with each Assessor appearing in 2 of them. It was shown, as well, that they were *not* isomorphic to the Octonions, as their O-trips' signing orders could never be made consistent with the "true" Octonions, as deviations always appear in at least 2 O-trip sequences, no matter how one transformed among labeling schemes. Raoul Cawagas has recently corroborated that all the primitive Sedenion zero-divisors indeed inhere in such algebraic objects, which he showed by his FINITAS algebraic software package correspond to a hitherto unidentified kind of loop-theory entity, which he calls a "quasi-octonion loop," or $\tilde{O}_L$. This has almost the same properties as the standard Octonion NAFIL (for Non-Associative Finite Invertible Loop), the 16-element maximal subloop (+ and – values of each axial unit), shorthanded $O_L$, of the Sedenions' 32-element loop. The signal difference between them is the former's failure to satisfy the Moufang identity. More, the Sedenions possess 7 copies of this non-Moufang loop, whose axis systems' imaginaries have indices identical to our "automorpheme" listings.

The Sedenions, as Cawagas showed, also possess 7 extra copies of the standard Octonion loop, each constructed by taking one of the 7 O-trips, then supplementing it with the 4 Sedenions whose indices are 8 and the XOR's of 8 with the O-trip. (Further subloops correspond to the smaller, quite familiar *groups*: the 35 Quaternion copies (7 O-trips + 28 S-trips); 15 4-cycles associated with the imaginaries; the 2-cycle of the plus and minus real units; and, the trivial 1-element additive group corresponding to 0.) Our septet of "8-balls" clearly appear as subloops within these, each appearing in 3 (1 per each O-trip containing its strut constant). This listing exhausts the groups and loops housed in the Sedenions, but not all the interesting algebras. The largest of the latter also have Octonion multiplication tables, with each 8-ball also appearing in 3 different ones; but as we'll soon see, they are far from identical to the creatures we've just looked at!

The new kind of algebra captures the breakdown of metric norms – and shows it as distinct from (and, in fact, prior to) the emergence of zero-divisors. These "lariats" are *adjustable* loops, hence metrically indeterminate, which helps explain their name; but their designation is also an acronym, for **L**ine **A**lgebras of **R**eal and **I**maginary **A**xis **T**ransforms, with the stress on the first word: for the entities whose behaviors are tracked in lariat multiplication tables are *oriented lines* comprising numerical axes, *not* length-fixed "units" residing at specifiable locations on same. The "stretchy" nature of zero divisor pairings, where the choice of points to multiply taken from each is arbitrary, is thereby captured . . . but without the zero-division property *itself* becoming manifest. Lariats, then, are like "sleeper cells": the zero-dividing "terrorism potential" of some of their elements may only become apparent when contexts shift, with no warnings given.

## 2. 8-Ball + Box-Kite = Switching Yard

It is a staple of silent-movie melodramas like "Perils of Pauline": the heroine is tied to the railroad tracks by a moustachioed villain in black; the hero arrives almost too late, but saves the day by throwing his weight against a large switching-yard lever arm which moves the rails (and thereby redirects the rushing train) just in the nick of time.

After contemplating the simplicity and structural richness of Box-Kites, it is hard to believe it had been assumed for a century that divisors of zero were pathological entities with nothing recommending their study. But now we've seen them, a child-like question arises: what happens if you change the edge-signs and multiply the "wrong" line elements together? If $(i_3 + i_{10})(i_6 - i_{15}) = 0$, what about $(i_3 + i_{10})(i_6 + i_{15})$? And what system of relations, if any, obtains between such "switched" products and our octahedral lattice?

As part of the answer, switch all the "–" signs to "+" in Box-Kite I's ABC Sail: we now orbit the vertices in near-group-theory style, in a 3-cycle, two per sail. Whether we pick the "/" diagonal as in the last paragraph, or the "\", as long as we're consistent, it doesn't matter: $kA \cdot kB = 2k^2 \cdot C$, $kB \cdot kC = 2k^2 \cdot A$, $kC \cdot kA = 2k^2 \cdot B$. Hence, if we set $k = \frac{1}{2}$, we almost have a group: $\frac{1}{2}(i_3 + i_{10}) \cdot \frac{1}{2}(i_6 + i_{15}) = \frac{1}{2}(i_5 + i_{12})$. The would-be group, though, is non-commutative, as signs on products flip when term orders reverse. But we don't have even a *non-Abelian* group, since recalibrating per above gives us $A^2 = B^2 = C^2 = -\frac{1}{2}$, which is clearly only half of what we need to get a Quaternion "working copy"!

To add to the confusion, picking proper "units" on such zero-divisor axes, so as to create isomorphs of the complex plane, requires a *different* calibration: the radius of the circle in the Argand diagram requires $k$ be set not to $\frac{1}{2}$, but the *square root* of $\frac{1}{2}$! This of course makes the case for reading a sign-switched "Sail" as a Quaternion copy even worse, since the 3 non-real axes can no longer be subgroups without being "tuned."

This is the motivation for the lariat construct: for if all calibrations be ignored, and only the placement and orientation of diagonal line elements be considered, we see that each Sail consists of 2 Quaternion-like lariats, with the caveat that instead of the real unit serving as *identity*, the (positively oriented) real axis *as a whole* does (call it a *lifeline*). This is readily shown true for all "sign-switched" Sails: the other 3 now have 2 "–" and 1 "+" edge-signs, but otherwise yield 2 sets of "Quizzical Quaternions" each.

Corresponding to William Rowan Hamilton's famous formula, we can write 2 "ijk"-style one-liners to capture this behavior for each Sail. Call the axis of reals "R," and use upper and lower case letters to stand for the (– –, ++) or "/" and (– +, + –) or "\" diagonals respectively: "A" now means $(i_3 + i_{10})$, while "a" stands for $(i_3 - i_{10})$. Show the *switched* edge-signs between the Sail's vertices, and here's what we get:

$$\begin{array}{ll}
A+B+C+: & A^2 = B^2 = C^2 = ABC = (-R) \\
& a^2 = b^2 = c^2 = abc = (-R) \\
\\
A-D+E-: & A^2 = d^2 = e^2 = Ade = (-R) \\
& a^2 = D^2 = E^2 = aDE = (-R) \\
\\
F-C-E+: & F^2 = c^2 = E^2 = FcE = (-R) \\
& f^2 = C^2 = e^2 = fCe = (-R) \\
\\
F+D-B-: & F^2 = D^2 = b^2 = FDb = (-R) \\
& f^2 = d^2 = B^2 = fdB = (-R)
\end{array}$$

It is curious to relate the "/" and "\" flips in the list just given to the correlated ordering of S-trip signings in each case. As mentioned at the end of the introductory section, each pair of Octonions in each Sail's O-trip can be replaced with their associated Sedenions, yielding an S-trip which circuits the same 3 Assessors. *Only with ABC*, however, whose edge-signs are all "–," are orientations maintained consistently: e.g., (3, 6, 5); (3, **15**, **12**); (**10**, 6, **12**); (**10**, **15**, 5) all cycle in left-right order. With the other Sails, 2 orientations are *always* reversed, with *only the S-trip whose Octonion is shared with ABC* displaying the same cycle order as its own O-trip. Striking enough to merit a name, as the "Trip-Sync Property," this fact should prove quite powerful in future simulation studies of emergent order, and "firefly tree"-type synchronization, especially among ensembles of Box-Kite-guided waves or particles – a theme we'll touch upon from other angles later on. For now, though, we return to lariats, to see what others the Sedenions might house.

If we include an 8-Ball with its proper Box-Kite, we find that the two diagonals at each end of a strut combine with the 8-Ball's "true copy" of the Quaternions, producing a lariat of "Mock Octonions." This is demonstrable by direct calculation, provided we set up good notation. Use "R" and upper vs. lower casing as per above; make "8" stand for the unit with that index, and "S" for the unit indexed by the strut constant; then, write "X" for the unit whose index is the "XOR" of 8 with the strut constant. If we next pick A and F as our strut pair, we can build a "Mock Octonion" multiplication table as follows: writing our symbols in the sequence (R, 8, X, S; F, a, f, A), we get the same table as if we'd associated these with the standard reals, and imaginaries indexed 1 – 7, respectively.

The same occurs if we substitute the other two struts: instead of "F, a, f, A" employ "E, b, e B" or "D, c, d, C" to the same effect . . .which means each pairing of Box-Kite with 8-Ball yields 3 isomorphic "Mock Octonion" lariats, as claimed earlier, for 21 in all. (As opposed to 2 x 4 x 7 = 56 "Quizzical Quaternion" lariats.) Here's the table:

**Generic Multiplication Table for "Mock Octonion" Lariats**

|   | *R* | *8* | *X* | *S* | *F* | *a* | *f* | *A* |
|---|---|---|---|---|---|---|---|---|
| *R* | +R | 8 | X | S | F | a | f | A |
| *8* | 8 | – R | S | – X | a | – F | – A | f |
| *X* | X | – S | – R | 8 | f | A | – F | – a |
| *S* | S | X | – 8 | – R | A | – f | a | – F |
| *F* | F | – a | – f | – A | – R | 8 | X | S |
| *a* | a | F | – A | f | – 8 | – R | – S | X |
| *f* | f | A | F | – a | – X | S | – R | – 8 |
| *A* | A | – f | a | F | – S | – X | 8 | – R |

The 8-Ball's Quaternion group is highlighted in the upper left quadrant; **F** and **A** stand in for the indices of the strut-pair represented by those letters on a given Box-Kite lattice. An isomorphic table results if **E** and **B** (and **e** and **b**), or **D** and **C** (and **d** and **c**), are substituted for **F** and **A** (and **f** and **a**) respectively. As with the Quizzicals, upper case, for a given Assessor, indicates the **(o + S)** or "/" diagonal; lower case, the **(o – S)** or "\".

# The Switching Yard

|   | *R* | *8* | *X* | *S* | *F* | *a* | *f* | *A* | *E* | *b* | *e* | *B* | *D* | *c* | *d* | *C* |
|---|---|---|---|---|---|---|---|---|---|---|---|---|---|---|---|---|
| *R* | +R | 8 | X | S | F | a | f | A | E | b | e | B | D | c | d | C |
| *8* | 8 | −R | S | −X | a | −F | −A | f | b | −E | −B | e | c | −D | −C | d |
| *X* | X | −S | −R | 8 | f | A | −F | −a | e | B | −E | −b | d | C | −D | −c |
| *S* | S | X | −8 | −R | A | −f | a | −F | B | −e | b | −E | C | −d | c | −D |
| *F* | F | −a | −f | −A | −R | 8 | X | S | −c | −D | 0 | 0 | b | E | 0 | 0 |
| *a* | a | F | −A | f | −8 | −R | −S | X | −D | c | 0 | 0 | E | −b | 0 | 0 |
| *f* | f | A | F | −a | −X | S | −R | −8 | 0 | 0 | −C | −d | 0 | 0 | B | e |
| *A* | A | −f | a | F | −S | −X | 8 | −R | 0 | 0 | −d | C | 0 | 0 | e | −B |
| *E* | E | −b | −e | −B | c | D | 0 | 0 | −R | 8 | X | S | −a | −F | 0 | 0 |
| *b* | b | E | −B | e | D | −c | 0 | 0 | −8 | −R | −S | X | −F | a | 0 | 0 |
| *e* | e | B | E | −b | 0 | 0 | C | d | −X | S | −R | −8 | 0 | 0 | −A | −f |
| *B* | B | −e | b | E | 0 | 0 | d | −C | −S | −X | 8 | −R | 0 | 0 | −f | A |
| *D* | D | −c | −d | −C | −b | −E | 0 | 0 | a | F | 0 | 0 | −R | 8 | X | S |
| *c* | c | D | −C | d | −E | b | 0 | 0 | F | −a | 0 | 0 | −8 | −R | −S | X |
| *d* | d | C | D | −c | 0 | 0 | −B | −e | 0 | 0 | A | f | −X | S | −R | −8 |
| *C* | C | −d | c | D | 0 | 0 | −e | B | 0 | 0 | f | −A | −S | −X | 8 | −R |

     Putting together the Box-Kite and 8-Ball associated with same strut constant, we get 1 of 7 isomorphic, 16-dimensional, "Switching Yards": the table given above can be used to instantiate them all, per the substitution rules given in the text. The 8-Ball "true Quaternion" subtable is shaded upper-left, and is contained in all 3 Mock Octonion lariats associated with a given strut constant. For easy parsing, the remaining ¾ of each such lariat are similarly color-coded. The zeros associated with the product-pairings along a Box-Kite's 12 edge-signs take up 48 of the table's 256 cells, and are painted pale yellow. (The simple task of eking out the octet of Quizzicals' identity-element-less "ijk-tables" from the remaining 8 x 6 white cells is left as an exercise for the curious.)

     In the above, zero-divisors have been sufficiently "tamed" to make contemplating their uses where most needed – in quantum theory. Some hints to this effect have been dropped in what's already been said. In the pages that remain, we'll pursue these further, plus other oddments which haven't been broached yet. It is implicit that "lariats" suggest a rather deeper answer to how Uncertainty works than just passing the buck to "probability waves"; for such an ascription to make sense, the vice of the "N Squares Rule" breakdown – shown long ago by Cayley's friend Kirkman to create a never-ending passalong to higher $2^n$-ions, if the product space Sedenions implicate is to be closed – must be turned into a virtue: which suggests some rethinking of what "renormalization" means, and why it is called for in the first place. And then there are some other things . . .

## 3. Switching-Yard Exotica: "Trigram Currents" and the "DJ Play-List Problem"

Say "+" on an edge-sign is seen as green, while "–" appears as red. If we "throw the switch" on all edge-signs at once, Sail ABC and Vent DEF reverse colors. But if we are viewing the square formed by any 2 struts from atop the third, changing "+" to "–" would appear like a 90º rotation (along all 3 struts simultaneously, in fact). What were called "tray-racks" in this series' first installment are the 4-cycles of zero-producing products we get when traversing the edges of any of the Box-Kite's 3 squares: B–C+E–D+; A–C+F–E+; A+D–F+C–. Toggling signs *shifts* sequences, but doesn't change them.

On the first effect, it's easy to see we get a 3-bit code uniquely telling us which Sail, in which switch-state, we're looking at. The patterns we get, in fact, are precisely those of the "coin oracle" method of creating I Ching readings. With 3 heads, or 2 tails and 1 head, we get a "yang" line: ABC, FDB, ADE, FCE give + + +, + – –, – + –, – – +, respectively, meaning we're in the land of "lariats," where products conceivably yield up indeterminate "measurement quanta" as allowed by Quaternionic or Octonionic representation schemes. But 3 tails, or 2 heads and 1 tail, yield a "yin" line instead: the "zero-division" condition, with the 4 Sails now bit-mapped to – – –, – + +, + – +, + + –.

If we then map "–" to 0, and "+" to 1, we get the 3-bit code, and then can start wondering about what sort of "breathing" would govern such yin-yang yoyo-ing. (And just as crucial would be grasping – perhaps via simulation studies of the "cellular automaton" variety – how 2 or more Box-Kite "particles" might interact, or which "hexagrams" might result from their exchanges under different conditions.) Given that the second paper in this series showed how the analogs of Box-Kites in higher $2^n$-ions can be decomposed or "folded" into the simpler forms found in the Sedenions, the study of such "trigram currents" might prove quite illuminating.

But now let's take a look at the second effect mentioned above: the "rotation" of tray-racks A–C+F–E+, A+D–F+C–, B–C+E–D+. Unlike the 2 types of 6-cycles (one, for A– B– C– only, called a "triple zigzag" due to its "/ \ / \ / \" alternation pattern; the other, obtaining for the other 3 Sails, dubbed a "trefoil" after the same-named knot), the alternating signs bring a tray-rack circuit back to the diagonal it started with in just 4 steps, so that *2* such circuits of the same quartet of Assessors (1 starting with a "/"diagonal, the other with a "\") are required to "cover all the bases."

But if we "throw the switch" to get lariats instead of zero-division sequencing, each tray-rack circuit generates 2 of the 4 terms comprising the $3^{rd}$ strut orthogonal to its plane in the Box-Kite octahedron – suggesting the cross-product familiar in vector calculus, but also having a unique higher-dimensional analog in the Octonions, which product is preserved in the Octonions' automorphism group, $G_2$ – which, since Moreno's work, we know governs zero-divisor structure, too, in the Sedenions as well as all higher $2^n$-ions.

Geoffrey Dixon has further shown how varieties of Octonion vector products can generate the Barnes-Wall and Leech lattices in 16- and 24- D, respectively. These not only play critical roles in "closest sphere-packing" and self-correcting error-code theories, but in M-Theory (the 16-D "4,320-packing" being the basis of the O(32) string). Even more compelling, Dominic Joyce's holonomy manifold with $G_2$ symmetry is currently seen as the most likely candidate for harboring the 7 "curled-up" gauge-dimensions of M-Theory. Making this seem even less a coincidence, the $G_2$ Joyce manifold is generated by a system of 7 x 7 equations which readily reduce to 35: viz., 7 O-trips, linked 1-to-1 with the 7 "8-Balls," with the triple-zigzag's substitution symmetry linking each such pairing uniquely to 3 S-trips, for 5 x 7 = 35 Quaternion copies in all, 1 per each Joyce generator?

To replace the "?" with a "." in that last sentence, we need to investigate what might be termed the "Dixon-Joyce circulatory system": we know we can get from the $G_2$ manifold to the Barnes-Wall lattice, hence M-Theory. We can imagine 35 equations each tracking one of the 35 Quaternion copies in the Sedenions. Which raises the "DJ Play-List Problem": how do we get from the "play list" of Barnes-Wall circuits Dixon's Octonion vector-products give us, to Sedenion orbits, and thence back down to Joyce's 35 choices? In the 8-D space of Octonions, there is an exact correspondence between the 240 "closest-packed" hyperspheres there and Octonion arithmetic, as Coxeter's classic study of "Integral Cayley numbers" made clear. But the link between closest-packing and number theory breaks down in higher dimensions; there are many different systems for navigating 16-D spaces, whose interrelationships are hardly as well-mapped as those in spaces of 8 dimensions or less.

The other 16-D-lattice-based string theory, "$E_8$ x $E_8$," is shown in M-Theory to be transformable into O(32), the latter comprising the objects at distance 2 on the same lattice whose unit-distance objects yield the famous anomaly-free "496" of the former: $E_8$ is just the "A, D, E Classification" code corresponding to the Cayley numbers, and twice 240 + 16-D Cartan subalgebra yields the "magic number." And this setup *does* have an easily identified (if not so well-known) hypernumber equivalence: to the so-called "M-algebra" posited by Charles Musès a quarter century ago. Just as the 4-space of reals plus Pauli spin-matrices is isomorphic to the Quaternions, Musès made use of the 2 x 2 identity matrix, with the usual imaginaries replacing the "1's," to commute between the 2 spaces (in effect, building the Clifford algebra $Cl_3$), then used the matrix to commute between the Octonions and an "extended Pauli" 8-space, which 2 spaces (by the same argument used in 4-D) are isomorphic. Their total 16-D arithmetic, then, "spans" $E_8$ x $E_8$ – albeit not as "nicely" as Octonions relate to integral Cayley numbers. What is of signal interest, however, is the fact that, real axis aside, the Sedenions are pure imaginaries, but Musès's M-algebra is a projection-operator-rich "half and half." Yet M-Theory transforms between the O(32) and $E_8$ x $E_8$ string theories. What clean transformations exist between Sedenion and "Switching Yard" arithmetics, and both Joyce and Dixon's domains?

One more 16-D space in need of interrelating arises when the icosahedral reflection group ("$H_3$" in "A, D, E" argot) is studied in singularity theory: members of Arnol'd's school showed a correspondence to a 16-D manifold with a 4-D boundary singularity – this latter so-called "open swallowtail" being a symplectic manifold in its own right . . . and first instance of an infinite family of symplectic invariants, per the Givental' theory of "triads." The link between symplectic and Quaternionic being self-evident, and the fact that Moreno's paper showed a 4-D "boundary" in Sedenion space, in the sense of "free of zero-divisors" (i.e., isomorphic to an "8-Ball"), one is led to wonder if there is some connection here – and if so, whether Gabrielov's triads relate, in turn, to the infinite sequence of $2^n$-ions (and hence, to some deep link between $G_2$ and $H_3$ in this context).

All of which suggests that the "pathology" of zero-divisors, far from marking the end in the road of useful extensions of the imaginary domain, may well mark the start of a new journey, to hitherto unvisited places. And, as lariats and whatever other "sleeper cell" structures hiding in higher dimensions can smuggle disguised zero-divisors into much lower-dimensioned algebras than suspected (with implications we can only begin to imagine for quantum foundations), so too we must allow for the possibility that much of what passed for Nature's architecture in 20[th] Century physics may prove to be removable scaffolding.

## 4. Post-Script on the "Trip-Sync Property," and Mathematical Beauty

Lab data displaying evidence of Chaos was, prior to that theory's emergence, typically tossed in the trash as so much "experimental error." Higher imaginaries containing zero divisors have likewise been dismissed, with the sort of superficial ugliness easily observed in the distribution of "Strut Table" S-trips – nicely mapped one per slot in ABC, it's true, but with 2 or 3 occurences of the same S-trip, with compensating absences of others, in the other 3 Sails – taken as endemic. But we've seen we can even contemplate a bonafide *representation theory* for zero divisors thanks to lariats, and find an astonishing degree of orderliness among S-trips thanks to the "Trip-Sync" property – the very property which underwrites the emergence of Sedenion zero-divisors in the first place. Beauty is not just in the mind of the beholder; it inheres in objective structures – but will only be perceived by those who open themselves to it without prejudice. Readers are invited to confirm the tabulation below with the evidence provided in the earlier listings of O-trips and S-trips, the better to rid themselves of any such aesthetic impediments.

**Synchronization Table By Box-Kite and Sail**
(Octonions in italics; trips in L-R cycle order shown in bold)

| BK# | ABC | ADE | FCE | FDB |
|---|---|---|---|---|
| I   | *(3 6 5)*    | *(3 4 7)*    | *(2 5 7)*    | *(2 4 6)*    |
|     | ( *3* 15 12) | ( *3* 13 14) | ( *2* 12 14) | ( *2* 13 15) |
|     | (10  6 12)   | (10  4 14)   | **(11  5 14)** | (11  4 15)   |
|     | (10 15  *5*) | (10 13  *7*) | (11 12  *7*) | **(11 13  *6*)** |
| II  | *(1 7 6)*    | *(1 4 5)*    | *(3 6 5)*    | *(3 4 7)*    |
|     | ( *1* 13 12) | ( *1* 14 15) | ( *3* 12 15) | ( *3* 14 13) |
|     | (11  7 12)   | (11  4 15)   | **( 9  6 15)** | ( 9  4 13)   |
|     | (11 13  *6*) | (11 14  *5*) | ( 9 12  *5*) | **( 9 14  *7*)** |
| III | *(2 5 7)*    | *(2 4 6)*    | *(1 7 6)*    | *(1 4 5)*    |
|     | ( *2* 14 12) | ( *2* 15 13) | ( *1* 12 13) | ( *1* 15 14) |
|     | ( 9  5 12)   | ( 9  4 13)   | **(10  7 13)** | (10  4 14)   |
|     | ( 9 14  *7*) | ( 9 15  *6*) | (10 12  *6*) | **(10 15  *5*)** |
| IV  | *(1 2 3)*    | *(1 7 6)*    | *(5 3 6)*    | *(5 7 2)*    |
|     | ( *1* 14 15) | ( *1* 11 10) | ( 5 15 10)   | ( 5 11 14)   |
|     | (13  *2* 15) | (13  *7* 10) | **( 9  *3* 10)** | ( 9  *7* 14) |
|     | (13 14  *3*) | (13 11  *6*) | ( 9 15  *6*) | **( 9 11  *2*)** |
| V   | *(2 4 6)*    | *(2 3 1)*    | *(7 6 1)*    | *(7 3 4)*    |
|     | ( *2*  9 11) | ( *2* 14 12) | ( 7 11 12)   | ( 7 14  9)   |
|     | (15  *4* 11) | (15  *3* 12) | **(10  *6* 12)** | (10  *3*  9) |
|     | (15  9  *6*) | (15 14  *1*) | (10 11  *1*) | **(10 14  *4*)** |
| VI  | *(3 4 7)*    | *(3 1 2)*    | *(5 7 2)*    | *(5 1 4)*    |
|     | ( *3* 10  9) | ( *3* 15 12) | ( 5  9 12)   | ( 5 15 10)   |
|     | (13  *4*  9) | (13  *1* 12) | **(11  *7* 12)** | (11  *1* 10) |
|     | (13 10  *7*) | (13 15  *2*) | (11  9  *2*) | **(11 15  *4*)** |
| VII | *(1 4 5)*    | *(1 2 3)*    | *(6 5 3)*    | *(6 2 4)*    |
|     | ( *1* 11 10) | ( *1* 13 12) | ( 6 10 12)   | ( 6 13 11)   |
|     | (14  *4* 10) | (14  *2* 12) | **( 9  *5* 12)** | ( 9  *2* 11) |
|     | (14 11  *5*) | (14 13  *3*) | ( 9 10  *3*) | **( 9 13  *4*)** |

"Trip-Sync" is not some "nice" feature obtaining in the Sedenions only. On the contrary, it seems to be a fundamental property of Box-Kites at whatever level, as one can readily verify in the 32-D world of the Pathions. The signal distinction between this arena and the 16-D Sedenions' is that the "emanation tables" therein – multiplication tables involving $2^{n-1} - 2$ entities, with one $2^n$-ion (whose index $s < 2^{n-1}$ is the strut constant) excluded from each – interconnect *multiple* Box-Kites, defined in the Pathions by 4-, not 3-bit, XORing. But for $s \leq 8$, each such 14-vertex entity decomposes into 7 Box-Kites, all showing "Trip-Sync" behavior. Take the O-trip associated with the ABC "triple-zigzag" in the Sedenion Box-Kite for the same $s$, then add '8' to the associated Sedenion in each Assessor to get the analogous Pathion. Using this Box-Kite as base-line, get the rest by successively keeping one strut-pair unchanged, but substituting for the other 2 pairs, in 1 of the 2 ways the 8 "extra" Assessors allow. For $s = 1$, the 14 Pathion Assessors' S and P indices (as opposed to Sedenions' **o ± S**) all have "inner XOR" of $16 + 1 = 17$, as follows:

 2,19  4,21  6,23  8,25  10,27  12,29  14,31  15,30  13,28  11,**26**  9,24  7,22  5,20  3,18

By the same "nested parentheses" logic used with Sedenions, successive pairings of outermost terms form struts, from (2, 19) and (3, 18) to (14, 31) and (15, 30). (See p. 14 of BK2 to find edge-signs given for the full "emanation table.") The 7 Box-Kites can be put in a table analogous to the "Strut Table" provided above – but, $s$ being the same for all, sequencing results from holding each of the 3 struts constant, while replacing the others in both ways that let ABC be an ASO S-trip, all 7 listed in ABC order:

|  | *A* | *B* | *C* | *D* | *E* | *F* |
|---|---|---|---|---|---|---|
| Base-Line Box-Kite: | 3 | 6 | 5 | 4 | 7 | 2 |
| AF kept constant: | 3 | 10 | 9 | 8 | 11 | 2 |
|  | 3 | 13 | 14 | 15 | 12 | 2 |
| CD kept (& moved to AF): | 5 | 12 | 9 | 8 | 13 | 4 |
|  | 5 | 14 | 11 | 10 | 15 | 4 |
| BE kept (& moved to AF): | 6 | 11 | 13 | 12 | 10 | 7 |
|  | 6 | 15 | 9 | 8 | 14 | 7 |

Use the P-trip lists in BK2 to confirm. When $s = 8$, the logic given fails, since no O-trip with an '8' can underwrite Sedenion ZD's. Here, however, what happens instead is that one still gets 7 Box-Kites, but the ABC Sail in each is uniquely associated with one of the 7 O-trips. For higher $s$, carry-bit "elbow room" runs out and "sand mandalas" with *only 3* Box-Kites result (see BK2), all 3 sharing strut ($s - 8$, 8). (Hence, using the *7* Sedenion Box-Kites as our "starter kit," we build 8 x 7 = *56*, plus 7 x 3 = *21*, for *84* Pathion Box-Kites all told: 1 per each Sedenion 8-Ball, Quizzical Quaternion and Mock Octonion lariat, and ZD diagonal, respectively.) For $s = 9$, for instance, we get this trio:

|  | *A* | *B* | *C* | *D* | *E* | *F* |
|---|---|---|---|---|---|---|
| Base-Line Box-Kite: | (**2**, 27) | (**8**, 17) | (**10**,19) | (3, 26) | (*1*, 24) | (11, 18) |
| "First Harmonic": | (4, 29) | (**8**, 17) | (12,21) | (5, 28) | (*1*, 24) | (13, 20) |
| "Second Harmonic": | (7, 30) | (**8**, 17) | (15,22) | (6, 31) | (*1*, 24) | (14, 23) |

Trip-Sync holds here too. Prove, qualify, or refute: it holds for *all* $2^n$-ions, $n > 3$.